\documentclass[10pt]{article}
\usepackage[utf8]{inputenc}
\usepackage[T1]{fontenc}
\usepackage[english]{babel}

\usepackage{lmodern} 
\usepackage[toc,page]{appendix}
\usepackage{fullpage} 
\usepackage{color}
\usepackage{amsmath} 
\usepackage{amsfonts} 
\usepackage{amssymb} 
\usepackage{amsthm} 
\usepackage[overload]{empheq} 
\usepackage{mathrsfs} 
\usepackage{pifont} 
\usepackage{verbatim} 

\usepackage{enumerate} 
\usepackage{stmaryrd}  
\DeclareMathAlphabet{\mathscrbf}{OMS}{mdugm}{b}{n}

\usepackage{boites} 
\usepackage{dsfont} 

\newcommand{\C} {\mathscr{C}} \newcommand{\R} {\mathbb{R}}
\newcommand{\Z} {\mathbb{Z}} \newcommand{\T} {{\mathbb T}} 
 
\newcommand{\uu} {\mathbf{u}} 
\newcommand{\xx} {\mathbf{x}} \newcommand{\xxi}{\boldsymbol{\xi}}
\newcommand{\na} {\nabla} 
\newcommand{\ffi}{\varphi} \newcommand{\ep} {\varepsilon}
\newcommand{\dd} {\mathrm{d}} 
\newcommand{\hh}{\mathbf{h}}

\newcommand{\conv}[2]{\operatorname*{\longrightarrow}_{#1 \rightarrow
    #2}} 
\newcommand{\convw}[2]{\operatorname*{\rightharpoonup}_{#1 \rightarrow
    #2}} \newcommand{\B}{\textnormal{B}}

\newcommand{\Ll}{\textnormal{L}_{\textnormal{loc}}}

\newcommand{\Hh}{\textnormal{H}} \renewcommand{\H}{\textnormal{H}}
\newcommand{\Hf}{\mathbf{H}} \newcommand{\Lf}{\mathbf{L}}
\newcommand{\W}{\textnormal{W}} \renewcommand{\P}{\textnormal{P}}
\newcommand{\D}{\mathscr{D}} 
 \newcommand{\Ss}{\textnormal{S}}
 \newcommand{\mm}{\mathbf{m}}
\newcommand{\uuin}{\uu_{\textnormal{in}}}
\newcommand{\vin}{v_{\textnormal{in}}}
\newcommand{\rhoin}{\rho_{\textnormal{in}}}
\newcommand{\fin}{f_{\textnormal{in}}}

\newcommand{\N}  {\mathbb{N}}

\newcommand{\yy} {\mathbf{y}}

\newcommand{\Pp} {\textnormal{P}} 
 \renewcommand{\P}{\mathbb{P}}
\newcommand{\vv} {\mathbf{v}} 
\newcommand{\oom} {\boldsymbol{\omega}} 
\newcommand{\aaa}{\mathbf{a}}

\renewcommand{\div}{\textnormal{div}}
\newcommand{\pa}{\partial}
\newcommand{\ds}{\displaystyle}

\newcommand{\beq}{\begin{equation}}
\newcommand{\eeq}{\end{equation}}
\newcommand{\ben}{\begin{eqnarray}}
\newcommand{\een}{\end{eqnarray}}
\newcommand{\beno}{\begin{eqnarray*}}
\newcommand{\eeno}{\end{eqnarray*}}
\newcommand{\beqno}{\begin{equation*}}
\newcommand{\eeqno}{\end{equation*}}
\newcommand{\be}{\begin{enumerate}}
\newcommand{\ee}{\end{enumerate}}
\newcommand{\bitem}{\begin{itemize}}
\newcommand{\eitem}{\end{itemize}}
\newcommand{\baln}{\begin{align}}
\newcommand{\ealn}{\end{align}}
\newcommand{\balno}{\begin{align*}}
\newcommand{\ealno}{\end{align*}}

\newtheorem{thm}{Theorem}[section]
\newtheorem{Lem}{Lemma}[section]
\newtheorem{rmk}{Remark}[section]

\newtheorem{prop}{Proposition}[section]

\newcommand{\Tt}{\T^3}
\newcommand{\Rt}{{\R^3}}

\newcommand{\Linf}{\textnormal{L}^{\infty}}

\newcommand{\boule}{\textnormal{\textsf{B}}}

\newcommand{\intt}[2]{\int_{#1}^{#2}}

\newcommand{\norm}[1]{\left\| {#1} \right\|}

\newcommand{\fiep}{\ffi_{\ep}} 
\newcommand{\gaep}{\gamma_{\ep}}
\newcommand{\fep}{f_{\ep}} 

\newcommand{\uep}{\uu_{\ep}}
\newcommand{\rhoep}{\rho_{\ep}}

\newcommand{\ptin}{\Pp_T^{\textnormal{in}}}
\newcommand{\jin}{J_{\textnormal{in}}}
\newcommand{\iin}{I_{\textnormal{in}}}
\newcommand{\ktin}{\textnormal{K}_T^{\textnormal{in}}}
\newcommand{\pt}{\Pp_T(\oom,h,\nu)}

\def\L{\textnormal{L}}
\usepackage{amsbsy}

\begin{document}

\begin{center}
{\bf{\large{Existence theory for the kinetic-fluid coupling when small droplets are 
treated as part of the fluid}}}
\end{center}
\bigskip

\bigskip

\begin{center}
S. BENJELLOUN\footnotemark[1], L. DESVILLETTES\footnotemark[1], AND A. MOUSSA\footnotemark[2]
\end{center}

\footnotetext[1]{\label{adr}CMLA, ENS Cachan \& CNRS, 61 av du Pr\'esident Wilson, 94235, Cachan, Cedex, France.}
\footnotetext[2]{UPMC Université Paris 06 \& CNRS, UMR 7598, LJLL, F-75005, Paris \& INRIA Paris-Rocquencourt, Équipe-Projet REO, BP 105, F-78153 Le Chesnay Cedex, France.}


\begin{abstract}
  We consider in this paper a spray constituted of an incompressible
  viscous gas and of small droplets which can breakup. This spray is
  modeled by the coupling (through a drag force term) of the
  incompressible Navier-Stokes equation and of the Vlasov-Boltzmann
  equation, together with a fragmentation kernel.  We first show at
  the formal level that if the droplets are very small after the
  breakup, then the solutions of this system converge towards the
  solution of a simplified system in which the small droplets produced
  by the breakup are treated as part of the fluid.  Then, existence of
  global weak solutions for this last system is shown to hold, thanks
  to the use of the DiPerna-Lions theory for singular transport
  equations.
\end{abstract}

\section{Introduction}

Sprays are complex flows which are constituted of an underlying gas in
which a population of droplets (or dust specks) are dispersed,
cf. \cite{these:ORourke}. There are various possibilities for modeling
such flows, depending in particular on the volume fraction of the
liquid phase (cf. \cite{eccomas} for example).
\par
We focus here on the case when the volume fraction occupied by the
droplets is small enough to be neglected in the equations (such sprays
are called thin sprays, cf. \cite{these:ORourke}), so that the
modeling of the liquid phase can be performed by the use of a pdf
(\emph{particles density function}) which solves Vlasov-Boltzmann
equation (cf. \cite{wil,code:kiva2}). Denoting $f:= f(t,\xx,
\xxi,r)\ge 0$ the number density of droplets of radius $r$ which at
time $t$ and point $\xx$ have velocity $\xxi$, the Vlasov equation
writes
\begin{equation}\label{vlasov}
\pa_t f + {\xxi} \cdot \nabla_{\xx} f + \nabla_{\xxi} \cdot (f \mathbf{\Gamma} ) = Q(f), 
\end{equation}
where $\mathbf{\Gamma}$ represents the acceleration felt by the
droplets (resulting from the drag force exerted by the gas), and $Q$
is an operator taking into account the complex phenomena happening at
the level of the droplets (collisions, coalescences, breakup).
\par
We also restrict ourselves to the case when the gas is incompressible
and viscous, which is for instance the usual framework when studying
the transport of sprays in the upper airways of the human lungs,
cf. \cite{these:Ayman}.
Accordingly, the gas is modeled by the incompressible Navier-Stokes
equation
\begin{align} &\nabla_\xx \cdot  \uu = 0, \\
  &\rho_g \left[\partial_t \uu + \nabla_{\xx} \cdot (\uu \otimes \uu)
  \right] + \nabla_{\xx} p - \mu \Delta_{\xx} \uu =
  \mathbf{F}_{\text{ret}} , \label{NS}
\end{align}
where $\rho_g$ is the constant density of the gas, $\uu:=\uu(t,\xx)
\in \R^3$ is its velocity, $\mu$ is its constant (dynamic) viscosity,
and $\mathbf{F}_{\text{ret}}$ is the retroaction of the drag force:
 \begin{equation}\label{ret}
   \mathbf{F}_{\text{ret}} (t,\xx) = - \int_0^{+\infty} \int_{\R^3} \frac43\,\rho_l \,r^3\,  f \, \mathbf{\Gamma} \, \dd \xxi \dd r.
 \end{equation}
 Finally, $\rho_l$ is the constant density of the liquid (so that the
 mass of the droplets of radius $r$ is $ \frac43\,\rho_l \,r^3$).
\par
For particles with small Reynolds numbers, the drag force is given by the simple formula (known as ``Stokes' law'')
 \begin{equation}\label{eq:trainee}
 \mathbf{\Gamma}(t,\xx,\xxi,r)  = -\frac92\frac{\mu}{\rho_l} \frac{ \xxi - \uu(t,\xx)}{r^2},
 \end{equation}
 that we shall systematically use in the sequel.
\par
The modeling of the breakup phenomena is important in the applications
and has led to various models appearing in the literature
(cf. \cite{code:kiva2}).  A typical form of the breakup kernel is
obtained when assuming that the droplets after breakup have the same
velocities as before breakup. The operator writes then
 \begin{eqnarray}\label{eqn:noyau frag}
 Q(f) (t,\xx,\xxi,r) = -\nu(\xxi,r)\, f(t,\xx,\xxi,r) + \int_{r^*>r} 
 b(r^*,r)\,
 \nu(\xxi,r^*) \, f(t,\xx,\xxi,r^*)\, dr^* ,
 \end{eqnarray}
where $\nu := \nu(\xxi,r)\ge 0$ is the fragmentation rate, and $b:=  b(r^*,r) \ge 0$ is related to the probability of ending up with droplets of radius $r$ out of the breakup of droplets of radius $r^*$.
 We finally obtain the following system 
 \begin{align} 
&\nabla_{\xx}\cdot \uu = 0, \label{sysor1}\\
&\partial_t f + \xxi \cdot \nabla_\xx f + \nabla_{\xxi} \cdot [ f\mathbf{\Gamma}] = Q(f), \label{sysor2}  \\
 \label{sysor3}
&\rho_g\,\Big[\partial_t\uu + \nabla_\xx \cdot (\uu \otimes \uu)\Big] + \nabla_\xx p - \mu \Delta_\xx \uu  = \mathbf{F}_{\text{ret}}.
\end{align}

Because of the dependence w.r.t. $r$ of $\mathbf{\Gamma}$ in eq. (\ref{eq:trainee}), 
we see that the drag force acting on the small droplets leads to very quick equilibration
of their velocity with the velocity of the gas. It is therefore natural to try
to write a set of equations replacing (\ref{ret}) -- (\ref{sysor3}), 
 in which the small droplets are considered as part of the gaseous phase.
\par
Indeed, in the context of the numerical simulation of eq. (\ref{vlasov}) thanks to a particle method, the small droplets which are produced
 because of the breakup can  lead to a high computational cost.
Once integrated (at the level of the continuous equations)  in the gaseous phase, they can be discretized along with the gas
thanks to a finite volume scheme, and their contribution to the computational cost remains in this way reasonable.
\bigskip

In order to perform a mathematical study of the system obtained by such an approximation, we consider the simplest
possible case, namely when the fragmentation rate $\nu$ takes the form $\nu(\xxi,r) := \tau^{-1}\mathds{1}_{r>r_2}$, for some constants $r_2>0$ and $\tau>0$, the latter being seen as a characteristic time of fragmentation.
We also assume that the aerosol is bidispersed: only two possible radii $r_1>r_2$ exist for the droplets, and the result of the breakup of particles of radius $r_1$ are particles of radius $r_2$. This implies that the density $f$ splits in the following way
\begin{align*}
f(t,\xx,\xxi,r) = f_1(t,\xx,\xxi)\,\delta_{r=r_1} + f_2(t,\xx,\xxi)\,\delta_{r=r_2},
\end{align*}
and that system  \eqref{ret} -- \eqref{sysor3} rewrites, after having
normalized all the constants (except $r_2$):
 \begin{align}[left= \empheqlbrace]
\label{sys:1_bis}&\nabla_{\xx}\cdot \uu = 0, \\
\label{sys:2_bis}&\partial_t f_1 + \xxi \cdot \nabla_\xx f_1 + \nabla_{\xxi} \cdot \left[ {f_1\,(\uu-\xxi)}\right] = -{f_1} ,   \\ 
\label{sys:3_bis}&\partial_t f_2 + \xxi \cdot \nabla_\xx f_2 + \nabla_{\xxi} \cdot \left[ \frac{f_2\,(\uu-\xxi)}{r_2^2}\right]= \frac{f_1}{r_2^3}
,   \\ 
\label{sys:4_bis}& \partial_t\uu + \nabla_\xx \cdot (\uu \otimes \uu) + \nabla_\xx p -  \Delta_\xx \uu  =
 - \int_{\R^3} (\uu-\xxi)\,f_1\,\dd \xxi -r_2 \int_{\R^3} (\uu-\xxi)\, f_2 \, \dd \xxi.
\end{align}
As we already explained, we are interested in the asymptotic regime when the particles resulting from breakup are becoming smaller and smaller. In our bidispersed model, this reduces to study the limit $r_2\to 0$.
\medskip

Denoting $\rho = r_2^3\displaystyle \int_{\R^3} f_2 \, \dd \xxi$, and noticing that
\begin{align*} \frac{\dd}{\dd t} \bigg\{ \int_{\T^3}\int_{\R^3} f_1\, \frac{|\xxi|^2}2\,  \dd \xxi \dd \xx &+
r_2^3\, \int_{\T^3}\int_{\R^3} f_2\, \frac{|\xxi|^2}2\,  \dd \xxi \dd \xx + \int_{\T^3} \frac{|\uu|^2}2\,\dd \xx \bigg\} \\
 &+ \int_{\T^3}\int_{\R^3} f_1\, {|\xxi - \uu|^2}\,  \dd \xxi \dd \xx\\
&+ \frac{1}{r_2^2}\, \int_{\T^3}\int_{\R^3} r_2^3f_2\, |\xxi - \uu|^2\,  \dd \xxi \dd \xx +  \int_{\T^3} |\nabla_{\xx} \uu |^2\,\dd \xx =0, 
\end{align*}
we see that (at the formal level) $r_2^3\, f_2(t, \xx, \xxi) \to \rho(t,x) \, \delta_{\xxi = \uu(t,\xx)}$ when $r_2 \to 0$.
\medskip

Integrating eq. (\ref{sys:3_bis}) against $r_2^3\, \dd \xxi$, we end up with
\begin{equation}\label{nn}
\partial_t \rho + \nabla_\xx \cdot [ \rho\, \uu ] = \int_{\R^3} f_1 \,\dd \xxi. 
\end{equation}
Then, integrating eq. (\ref{sys:3_bis}) against $r_2^3\, \xxi\, \dd \xxi$ and adding the result with eq. (\ref{sys:4_bis}), 
we obtain 
$$ \partial_t ( (1 + \rho)\, \uu) + \nabla_\xx \cdot ( (1+\rho)\, \uu \otimes \uu) + \nabla_\xx p - \Delta_\xx \uu  = -
 \int_{\R^3}(\uu - \xxi)\, f_1\,\dd\xxi + \int_{\R^3} f_1\,\xxi\, \dd\xxi . $$
Combining this last equation with eq. (\ref{nn}) and replacing the notation $f_1$ by $f$, we write down the
 system that we wish to study (we close it with periodic boundary conditions for the mathematical proof of existence):
 \begin{align}[left= \empheqlbrace]
&\nabla_{\xx}\, \uu = 0, \label{sys_limite:1_bis}\\
&\partial_t \rho + \nabla_\xx \cdot [ \rho \, \uu ] = \int_{\R^3} f\, \dd \xxi, \label{sys_limite:3_bis}\\
&\partial_t f + \xxi \cdot \nabla_\xx f + \nabla_{\xxi} \cdot [ (\uu-\xxi)\, f] = -f,  \label{sys_limite:2_bis} \\ 
& (1+\rho)\Big[\partial_t\uu + \nabla_\xx \cdot (\uu \otimes \uu)\Big] + \nabla_\xx p - \Delta_\xx \uu  = 2\,
 \int_{\R^3}(\xxi-\uu)\, f\,\dd\xxi,\label{sys_limite:4_bis}
\end{align}
where $\rho:= \rho(t, \xx) \ge 0$, $\uu := \uu(t, \xx) \in \R^3$, $p:= p(t,\xx)\ge 0$, 
$f := f(t,\xx, \xxi) \ge 0$, and $t\ge 0$, $\xx\in \Tt$, $\xxi \in \R^3$. Let us recall
 that $\rho$ represents here the ``added density'' resulting from the very small particles. This is the reason why, though we normalized the fluid density in \eqref{sys:4_bis} by $\rho_g =1$, we have in \eqref{sys_limite:4_bis} the term $(1+\rho)$ in front of the convective part of the fluid equation. For a more detailed version of the previous computation, see \cite{these:Saad}.

\bigskip 

Our goal is to study the existence theory for this system, completed with the following initial data:
\begin{align}
f(0,\xx,\xxi) &= \fin(\xx,\xxi) \geq 0,\quad \xx \in \Tt,\quad \xxi \in \Rt, \label{ci:f}\\
\rho(0,\xx) &= \rhoin(\xx) \geq 0,\quad \xx \in \Tt, \label{ci:rho}\\
\uu(0,\xx) &= \uuin(\xx)\in \Tt,\quad \xx\in \Tt. \label{ci:u}
\end{align}

Before stating our result, let us introduce a few notations. If there is no ambiguity on the considered time interval, we simply denote by $\L^p_{t}(\L^q_\xx)$ and $\L^p_t(\L^q_{\xx,\xxi})$ the spaces $\L^p\big([0,T];\L^q(\T^3)\big)$ and $\L^p\big([0,T];\L^q(\T^3\times\R^3)\big)$, for any pair of exponents $(p,q)\in[1,\infty]^2$. In particular
when $p=q$ we simply use the notation $\L^p$. We adopt the same convention for Sobolev spaces $\W^{1,p}_{\xx}$ and $\H^m_{\xx}$ ($m\in\N$) and denote by $\H^{-m}_{\xx}$ (or equivalently $\H^{-m}(\T^3)$) the dual of the latter.
%
%
 When the subscript `` div\,'' is added to any space, 
the corresponding subspace are the divergence-free (in $\xx$ only) elements of the ambient space. We denote by $\mathbb{P}$ the Leray projector onto divergence-free vector fields
\begin{align*}
\mathbb{P} = \textbf{Id} - \nabla_{\xx} \Delta^{-1}_{\xx}\div_{\xx},
\end{align*}
easily defined on $\mathscrbf{D}([0,T]\times\T^3)$ thanks to Fourier analysis, extended by duality to $\mathscrbf{D}'([0,T]\times\T^3)$. 

\vspace{2mm}

Finally, if $h$ is a scalar function defined on $ \R_+ \times\T^3 \times\R^3$, we define the following moments for $h: = h(t,\xx,\xxi) $ (and for $\alpha \in \R$)
\begin{align*}
m_\alpha h(t,\xx):=\int_{\Rt} |\xxi|^\alpha\, h(t,\xx,\xxi)\,\dd \xxi,
\quad M_\alpha h (t):= \int_{\T^3} m_\alpha(h)(t,\xx)\,\dd \xx, \quad\mm_1 h(t,\xx):=\int_{\Rt}\xxi h(t,\xx,\xxi)\, \dd \xxi.
\end{align*}
\smallskip



Thanks to those notations, we are able to write down our main
Theorem:

\begin{thm}
\label{thm:sys_2} 
Let $T>0$ and assume that  
\begin{equation}\label{ci:hypo}
\fin \in \textnormal{L}^{\infty},  (1+|\xxi|^2)\fin \in\L^1, 
 \rhoin \in \textnormal{L}^{\infty},
 \uuin\in \Lf_{\div}^2.
\end{equation}
Then the system \eqref{sys_limite:1_bis} -- \eqref{ci:u}
admits a global weak solution $(\rho \ge 0, \, f \ge 0, \,\uu \in \R^3)$ such that
\begin{equation}\label{esp}
\rho \in \L^{\infty}([0,T],\L^{5/3}(\Tt)),
f \in \L^{\infty}([0,T] \times \Tt \times \Rt),
\uu \in \L^{2}([0,T],\Hh^{1}_{\textnormal{div}}(\Tt)) \cap \L^{\infty}([0,T],\L^{2}(\Tt)).
\end{equation}
Moreover, the triplet $(\rho, \, f, \,\uu)$ satisfies the following energy estimate
\begin{align}
\frac{1}{2}\left\{M_2 f (t) + \| \sqrt{1+\rho(t)}\, \uu(t)\|_{\textnormal{L}^2(\T^3)}^2\right\} &+ \int_0^t\|\nabla_\xx \uu(s)\|_{\textnormal{L}^2(\T^3)}^2 \dd s + \frac{3}{2}\int_0^t\int_{\T^3 \times \R^3}|\uu-\xxi|^2 f \dd \xxi \dd\xx \dd s  \nonumber \\
& \leq \frac{1}{2}\left\{ M_2 f_{\text{\bf in}} + \| \sqrt{1+\rho_{\text{\bf in}}}\, \uu_{\text{\bf in}}\|_{\textnormal{L}^2(\T^3)}^2 \right \},
\end{align}
and the bound
\begin{align*}
\| f \|_{\Linf \left( [0,\,T] \times \Tt \times \Rt \right) } \leq e^{2T} \| f_{\text{\bf in}} \|_{\Linf\left( \Tt \times \Rt \right)}.
\end{align*} \end{thm}

\begin{rmk} \label{we}
We explain here the meaning of `` weak solutions '' in the above Theorem. For $p\in[1,\infty]$, when $\uu\in\Lf^p$, $f\in\L^\infty$ and $\rho\in \L^{p'}$,
 equations \eqref{sys_limite:1_bis}--\eqref{sys_limite:3_bis} have a clear meaning in the distributional sense. However the meaning of \eqref{sys_limite:4_bis} is not completely obvious, and we rewrite this equation, using both \eqref{sys_limite:4_bis} and \eqref{sys_limite:2_bis}, as
\begin{align}
\label{sys_limite:5_bis}\P\left\{\partial_t[(1+\rho)\uu] + \div_\xx \Big[(1+\rho) \uu \otimes \uu\Big]\right\} -\Delta_\xx \uu  = \P \left\{ 2\mm_1 f-\uu m_0 f \right\}.
\end{align}
For the initial condition, we adopt the following general definition : consider $\vin\in\D'(\T^3)$ and $w$ in $\Ll^1([0,T];\H^{-m}(\T^3))$ for some integer $m\in\N$. A
 distribution $v\in\D'\big(]-\infty,T[\times\T^3\big)$ is a solution of the Cauchy problem
$$ \partial_t v  = w, \qquad v(0) = \vin, $$
if the support of $v$ is included in $\R_+$ and if
\begin{align*}
\partial_t v  &= \widetilde{w} + \delta_0 \otimes \vin \text{ in }\D'\big(]-\infty,T[\times\T^3\big),
\end{align*}
where $\widetilde{w}$ is the extension of $w$ by $0$ on $\R_-$.
\end{rmk}
\medskip

The study of the existence of solutions to coupled (through drag force interaction) fluid-kinetic equations
is now a well-established subject.
\par
 As far as viscous equations are concerned for the modeling of the gaseous phase, we would like to quote the works on the Vlasov-Stokes equations in  \cite{art:Hamdache_stokes}, on the Vlasov/incompressible Navier-Stokes equations in \cite{art:Anoshchenko_VNS,article:Ayman,art:chengyu}, on the Vlasov-Fokker-Planck/incompressible Euler equations in \cite{article:caridumou}
, on the Vlasov-Fokker-Planck/incompressible Navier-Stokes equations in  \cite{art:Vasseur_comp_exist}, and on the Vlasov-Fokker-Planck/compressible Navier-Stokes equations in \cite{art:Vasseur_comp_limit}. 
\par
Without any viscosity or Fokker-Planck damping, the study is more difficult, and only smooth local solutions are known to exist, in the compressible setting (cf. \cite{art:Clne} and \cite{article:Jlien}). 
\bigskip

Our system \eqref{sys_limite:1_bis} -- \eqref{sys_limite:4_bis} is reminiscent of
the Vlasov/incompressible Navier-Stokes equations with a variable density.
It is known that variable densities in the incompressible Navier-Stokes 
equations lead to extra difficulties (w.r.t. constant densities), cf.  \cite{book:Boyer,book:Plions_incomp}. Here extra difficulties (w.r.t. \cite{article:Ayman}) appear in the final passage to the limit (that is, in the stability result for solutions satisfying the natural {\sl{a priori}} estimates
of the problem). 
They are linked with nonlinearities which are specific of our model, and which necessitate a
careful treatment using refined versions of Lemmas presented in \cite{art:diperna_lions}. 
\par
The equations coming out of the theory of sprays are known to be difficult to approximate in a good way (that is, in a way in which the {\sl{a priori}} estimates satisfied by the 
equations themselves, such as the energy estimate, are also satisfied or well approximated by the approximating equations). This difficulty appears in this paper, and the approximating scheme which is used is rather complicated and necessitates successive steps, which are 
reminiscent of those used in \cite{article:Ayman}. 
\bigskip

Our paper is structured as follows:
 Section 2 is devoted to the Proof of existence  of solutions to an approximated
version of this system. The approximation is  
removed in Section 3 and Theorem \ref{thm:sys_2}
is proven there. Finally, we present in a short Appendix some autonomous results 
which are used in the Proof of Theorem \ref{thm:sys_2}.

\section{Existence for a regularized system}

In this section, we begin the Proof of Theorem \ref{thm:sys_2} by introducing a regularized system.

\subsection{Definition of the regularized system}

In order to do so,
we introduce a mollifier $\ffi_{\ep}\in\mathscr{C}^\infty(\T^3)$ for all $\ep>0$.
We also introduce a truncation in this way: we define $\gamma_{\ep}$ a nonnegative element of $\D(\R^3)$,
 whose support lies in $\boule(0,2/\ep)$, bounded by $1$, and equal to $1$ on the ball $\boule(0,1/\ep)$. 
\medskip

Our regularized system writes
\begin{align}[left={\empheqlbrace}]
\label{eq:sysreg1} &\div_{\xx}\, \uep = 0, \\
\label{eq:sysreg2} &\partial_t\rhoep +(\uep\star\ffi_{\ep})\cdot \nabla_\xx\,\rhoep =m_0(f_\ep \gamma_\ep),\\
\label{eq:sysreg3} &\partial_t \fep + \xxi \cdot \nabla_\xx \fep +\nabla_{\xxi} \cdot [(\uep\star\ffi_{\ep}-\xxi)\fep]=-\fep,\\
\label{eq:sysreg4} &\mathbb{P}\Big\{(1+\rho_\ep)\big[\partial_t\uu_\ep + (\uu_\ep\star\ffi_\ep)\cdot \nabla_\xx \uu_\ep\big]\Big\} -\Delta_\xx \uu_\ep = 2\P\Big\{\mm_1(f_\ep\gamma_\ep)-\uu_\ep m_0(f_\ep\gamma_\ep)\Big\},
\end{align}
with initial conditions  $\fin^\ep\in\D(\T^3\times\R^3)$, $\rhoin^\ep\in\mathscr{C}^\infty(\T^3)$ and $\uuin^\ep\in\mathscrbf{C}^\infty_\div(\T^3)$ approximating respectively $\fin$ in all $\Ll^p(\T^3\times\R^3)$ ($p<\infty$) and in $\L^1(\T^3\times\R^3, (1+|\xxi|^2)\,\dd\xx\,\dd\xxi)$, $\rhoin$ in all $\L^p(\T^3)$ ($p<\infty$), and $\uuin$ in $\Lf^2(\T^3)$.
\medskip

Next subsection is devoted to the Proof of existence of solutions to this approximated system (that is, for a given 
parameter $\ep>0$).

\subsection{Existence of solutions for the regularized system} \label{sub22}

In this Subsection, the parameter $\ep$ is fixed, and we drop out the corresponding indices in order to make
 the formulas more readable.
\medskip

 We fix a triplet of initial data $(\oom,h,\nu)\in \mathscrbf{C}^\infty_\div\times ( \mathscr{C}^1\cap\L^\infty ) \times \mathscr{C}^1 $ (the two last ones being nonnegative). The considered system (in which we do not explicitly write the initial data) is then the following:
\begin{align}[left={\empheqlbrace}]
\label{eq:vnsadiv_} &\div_\xx\, \uu = 0,\\
\label{eq:vnsarho_}&\partial_t\rho +(\uu\star\ffi)\cdot \nabla_\xx\,\rho = m_0(f\gamma),\\
\label{eq:vnsaf_}&\partial_t f + \xxi \cdot \nabla_\xx f +\nabla_{\xxi} \cdot [(\uu\star\ffi-\xxi) f]=- f,\\
\label{eq:vnsau_} &\mathbb{P}\Big\{(1+\rho)\big[\partial_t\uu + (\uu\star\ffi)\cdot \nabla_\xx \uu\big]\Big\} -\Delta_\xx \uu = 2\P\Big\{\mm_1(f\gamma)-\uu m_0(f\gamma)\Big\}.
\end{align}

We are able to prove the
\medskip

\begin{prop} \label{regul}
For any triplet of initial data $(\oom,h,\nu)\in \mathscrbf{C}^\infty_\div\times ( \mathscr{C}^1\cap\L^\infty ) \times \mathscr{C}^1 $ (the two last ones being nonnegative), there exists a weak solution (in the sense of Remark \ref{we}) 
 to the system (\ref{eq:vnsadiv_}) -- 
(\ref{eq:vnsau_}).
\end{prop}

{\bf{Proof of Proposition \ref{regul}}}:

We shall show existence of a solution to the nonlinear system \eqref{eq:vnsadiv_} -- \eqref{eq:vnsau_} by applying Schauder's fixed point Theorem to the following mapping:
\begin{align*} 
\Ss:\mathscr{C}^0_t(\Lf^2_\div)&\longrightarrow \mathscr{C}^0_t(\Lf^2_\div)\\
\uu &\longmapsto \uu^\diamond,
\end{align*}
where $\uu^\diamond$ will be built thanks to the following steps :

\begin{enumerate}[(i)]
\item We first consider $f^\diamond\in\mathscr{C}^1$,  unique classical solution of 
\begin{align}
\label{eq:fixpointf}\partial_t f^\diamond + \xxi\cdot\nabla_\xx f^\diamond + \nabla_{\xxi}\cdot[(\uu\star\ffi -\xxi)f^\diamond] = -f^\diamond,
\end{align}
with $h$ as initial datum,
\item then $\rho^\diamond\in\mathscr{C}^1$ is defined as the unique classical solution of
\begin{align}
\label{eq:fixpointrho}\partial_t \rho^\diamond + (\uu\star\ffi)\cdot \nabla_{\xx}\rho^\diamond = m_0(f^\diamond\gamma),
\end{align} 
with  $\nu$ as initial datum,
\item and finally  $\uu^\diamond\in\mathscr{C}^0_t(\Lf^2_\div)$ is built as the unique divergence-free weak solution of 
\begin{align}
\hspace{-2.3cm}\label{eq:fixpointu}&\mathbb{P}\Bigg\{(1+\rho^\diamond)\Big[\partial_t\uu^\diamond + (\uu\star\ffi)\cdot \nabla_\xx \uu^\diamond\Big]\Bigg\} -\Delta_\xx \uu^\diamond = 2\P\Big\{\mm_1(f^\diamond\gamma)-\uu^\diamond m_0(f^\diamond\gamma)\Big\},
\end{align}
with  $\oom$ as initial datum. 
\end{enumerate}
In the sequel $\Pp_T(\oom,h,\nu)$ will denote  any nonnegative function (which also may depend on $\ep$, but as explained before we omit the corresponding index here) that splits into a finite sum of positively homogeneous functions of strictly positive degree  of one of the following terms : $\|\oom\|_{\Hf^1}$,$\|h\|_{\infty}$ or $\|\nu\|_{\infty}$. $\pt$ may change from one line to another,
 but will always have the structure that we just described.
\smallskip

In the next paragraph, we show that $f^\diamond$, $u^\diamond$ and $\rho^\diamond$ are well defined.
 
\subsubsection{Existence and uniqueness of $f^\diamond$, $\rho^\diamond$ and $\uu^\diamond$}\label{subsec:lipfrho}

We know that $\uu\star\ffi\in\mathscr{C}^0_{t}(\mathscrbf{C}^1_\xx)$ and $\nabla_{\xx} (\uu \star \ffi) \in\Lf^\infty$, so that  the characteristic curves of \eqref{eq:fixpointf} and \eqref{eq:fixpointrho}
 are globally well-defined.  We hence easily obtain the existence and uniqueness of two nonnegative
 functions $f^\diamond$ and $\rho^\diamond$ (classically) solving
\eqref{eq:fixpointf} and \eqref{eq:fixpointrho}. Thanks to the maximum principle, we also know that
\begin{align*}
\|f^\diamond\|_{\infty}\leq  e^{2T}\|h\|_{\infty},
\end{align*}
from which we get
\begin{align}
\label{eq:estimoment} \|m_0(f^\diamond\gamma)\|_{\infty} + \|\mm_1(f^\diamond\gamma)\|_{\infty} \leq C_\gamma e^{2T}\|h\|_\infty \leq \pt,
\end{align}
and hence
\begin{align}
\label{eq:estirho}\|\rho^\diamond\|_{\infty} &\leq \pt.
\end{align}

\vspace{2mm}

We now consider $f^\diamond$ and $\rho^\diamond$ as given. The existence and uniqueness of $\uu^\diamond \in\Hf^1_{t,\xx}\subset\mathscr{C}_t^0(\Lf^2_{\xx})$,  divergence-free weak solution of \eqref{eq:fixpointu} (with $\oom$ for initial data) may be obtained thanks to the usual Galerkin method. Since this type of construction is rather standard, we won't detail it here, and refer to \cite{these:Saad} for a precise treatment of this procedure in our system, or for 
instance to \cite{book:anto}, Chapter~3, for a more generic approach.
 \medskip
 
 In next subsection, we obtain natural estimates for the quantity $\uu^\diamond$. 

\subsubsection{Estimates for $\uu^\diamond$}\label{subsec:lipfrho_estimates}

 Taking $\uu^\diamond$ as test function in \eqref{eq:fixpointu} and using \eqref{eq:fixpointrho},
 we obtain the following estimate (since $\uu\star\ffi$ is divergence free and $f^\diamond$ is nonnegative)
\begin{align}
\nonumber\frac{1}{2} \,\left\{\int_{\T^3}(1+\rho^\diamond(t,\xx))|\uu^\diamond(t,\xx)|^2\dd \xx\right\} &+ \int_0^t\int_{\T^3}|\nabla_\xx \uu^\diamond(s,\xx)|^2\dd \xx\,\dd s \\
\nonumber&\leq-\frac{3}{2}\, \int_0^t\int_{\T^3}|\uu^\diamond(s,\xx)|^2 m_0(f^\diamond\gamma)(s,\xx)\,\dd \xx\,\dd s\\
\nonumber&+2\int_0^t\int_{\T^3}\uu^\diamond(s,\xx)\cdot \mm_1(f^\diamond\gamma)(s,\xx)\,\dd \xx\,\dd s\\
\label{ineq:nrjep}&+\frac{1}{2} \,\left\{\int_{\T^3}(1+\nu(\xx))|\oom(\xx)|^2\dd \xx\right\},
\end{align}
so that using $\rho^\diamond \ge 0$ and Gr\"onwall's lemma together with
 estimates  \eqref{eq:estimoment}--\eqref{eq:estirho}, we get 
\begin{align}
\label{eq:estimation1}
\|\uu^\diamond\|_{\L^\infty_t(\Lf^2_\xx)} + \|\uu^\diamond\|_{\L^2_t(\Hf^1_\xx)} \leq \pt\exp(\pt).
\end{align}
Taking then $\partial_t \uu^\diamond$ as a test function in \eqref{eq:fixpointu},
 we get (using again $\rho^\diamond \ge 0$ and $f^\diamond \ge 0$)
\begin{align*}
\|\partial_t \uu^\diamond\|_2^2 + \frac{1}{2}\,\|\nabla_{\xx}\uu^\diamond(t)\|_2^2 &\leq \|(1+\rho^\diamond) (\uu \star\ffi)\|_\infty \,\|\nabla_{\xx} \uu^\diamond\|_2 \, \|\partial_t \uu^\diamond\|_2 + 2\,\|\mm_1(f^\diamond \gamma)\|_2 \, \|\partial_t \uu^\diamond\|_2\\
& + 2\,\|m_0(f^\diamond\gamma)\|_\infty\, \|\uu^\diamond\|_2\, \|\partial_t\uu^\diamond\|_2  + \frac{1}{2}\, \|\nabla_{\xx}\oom\|_2^2,
\end{align*}
so that using Young's inequality, estimate \eqref{eq:estimation1} above,
 and estimates  \eqref{eq:estimoment}--\eqref{eq:estirho}, we end up with 
\begin{align}
\label{eq:estimation2}\|\partial_t \uu^\diamond\|_{2}  \leq \pt\,\exp(\pt)\,
 \Big[\|\uu\|_{\L^\infty_t(\Lf^2_\xx)}+1\Big],
\end{align}
thanks to the elementary convolution inequality
 $\|\uu\star\ffi\|_\infty \leq \|\uu\|_{\L^\infty_t(\Lf^2_{\xx})} \|\ffi\|_2$.
\medskip

Given any triplet $(\oom,h,\nu)\in \mathscrbf{C}^\infty_\div\times \mathscr{C}^1\cap\L^\infty\times\mathscr{C}^1$ of initial data, we can then define the mapping
\begin{align*} 
\Ss:\mathscr{C}^0_t(\Lf^2_\div)&\longrightarrow \mathscr{C}^0_t(\Lf^2_\div)\\
\uu &\longmapsto \uu^\diamond .
\end{align*}
\smallskip

In next paragraph, we use Schauder's fixed point Theorem for this mapping.

\subsubsection{Application of Schauder's Theorem}

According to the estimates of the previous paragraph, we see that
 for any $\uu\in\mathscr{C}^0_t(\Lf^2_{\div})$,
\begin{align*}
\|\Ss(\uu)\|_{\L^\infty_t(\Lf^2_{\xx})}\leq  \pt\exp(\pt):=R.
\end{align*}
Hence, $\Ss$ sends $\B_R$ (the closed ball of radius $R$ of the normed space $\mathscr{C}^0_t(\Lf^2_{\div})$) to itself. 
\medskip

Thanks to estimates \eqref{eq:estimation1} and \eqref{eq:estimation2}, a bounded subset of $\mathscr{C}^0_t(\Lf^2_\div)$ is sent by $\Ss$ to a bounded subset of $\Hf^1_{t,\xx}$, and Lemma \ref{annex:aubin} of the Appendix allows us to conclude that $\Ss(\B_R)$ is relatively compact in $\mathscr{C}^0_t(\Lf^2_\div)$.
\medskip

It remains to study the continuity of $\Ss$.
Notice first that for any $\uu\in\mathscr{C}^0_t(\Lf^2_{\div})$, any weak bounded solution of \eqref{eq:fixpointf} or \eqref{eq:fixpointrho} is a renormalized solution in the DiPerna-Lions sense introduced in \cite{art:diperna_lions}: the absorption coefficient ($=1$) of the first equation and the right-hand side of second one are bounded and one easily checks that the vector fields $(t,\xx,\xxi)\mapsto (\xxi,\uu\star\ffi(t,\xx)-\xxi)$ and $(t,\xx)\mapsto \uu \star \ffi(t,\xx)$ satisfy all the desired assumptions given in \cite{art:diperna_lions}. 
Now if $(\uu_n)_{n}$ converges to $\uu$ in $\mathscr{C}^0_t(\Lf^2_\div)$, we know from the previous step that $(\Ss(\uu_n))_{n\in\N}$ has a converging subsequence. Showing that the whole sequence converges to $\Ss(\uu)$ reduces hence to prove that it has only one accumulation point, namely $\Ss(\uu)$. Assume therefore that $\vv\in\mathscr{C}^0_t(\Lf^2_{\div})$ is an accumulation point of $(\Ss(\uu_n))_{n\in\N}$ and denote by $\sigma$ the corresponding extraction. Since $(\uu_{\sigma(n)})_n$ is still bounded, using \eqref{eq:estimation1}--\eqref{eq:estimation2} and adding a subsequence if necessary,
 we may assume that $\vv$ is also the limit of $(\Ss(\uu_n))_{n\in\N}$ in $\Hf^1_{t,\xx}$, for the weak topology. On the other hand, since $(\uu_{\sigma(n)})_n$ still  converges to $\uu$ in $\mathscr{C}^0_t(\Lf^2_{\div})$, DiPerna-Lions stability result of \cite{art:diperna_lions} ensures that $(f_{\sigma(n)})_n$ and $(\rho_{\sigma(n)})_n$ both converge in $\L^p_{\textnormal{loc}}$
 (for all $p<\infty$) to $f$ and $\rho$, the corresponding (unique) weak and bounded solutions of the associated equations (with initial data $h$ and $\nu$ and vector field defined by $\uu$). At this point, we have enough (strong) convergences to ensure that all the nonlinear terms of the fluid equation indeed
 converge to the expected limit : $\vv\in\Hf^1_{t,\xx}$ is hence a weak solution of \eqref{eq:fixpointu} with $\oom$ as initial datum. We already mentioned the uniqueness of solutions for this equation, so that we  eventually get $\vv = \uu^\diamond = \Ss(\uu)$. 
\medskip

$\Ss$ is hence a continuous map from a closed convex nonempty
 set of a normed space to itself, and it has a relatively compact range.
Thanks to Schauder's fixed point Theorem (see \cite{giltru} for instance),  $\Ss$ has a fixed point which is a 
 solution to \eqref{eq:vnsadiv_} -- \eqref{eq:vnsau_} with initial data $(\oom,h,\nu)$.
This concludes the Proof of Proposition \ref{regul}. $\square$
\medskip

Reintroducing the parameter $\ep>0$, we have obtained the existence (for any $T>0$, $\ep>0$) of
 a divergence-free vector field $\uep \in \Hh^1_{t, \xx}$, and two nonnegative functions $\fep \in \mathscr{C}^1_{t,\xx,\xxi}$ and $\rhoep \in \mathscr{C}^1_{t,\xx}$ solutions of system \eqref{eq:sysreg1} -- \eqref{eq:sysreg4}. 
\smallskip

Note that since $\uu_\ep\star\ffi_\ep\in\L^\infty_{t,\xx}$ and $\fin^\ep$ is compactly supported, one can
show without difficulty that $f^\ep(t,\cdot,\cdot)$ remains compactly supported for a.e. $t$ (with a support depending on $\ep$).
\medskip

In next section, we pass to the limit when $\ep \to 0$ in the functions $\uu_\ep$, $f_\ep$ and $\rho_\ep$, and show that their limit is a (weak) solution of \eqref{sys_limite:1_bis} -- \eqref{sys_limite:4_bis}.

\section{Passage to the limit and Proof of the main Theorem}

We keep on proving Theorem \ref{thm:sys_2} in this Section. We start with a local result of existence:

\begin{prop}
\label{prop:local} Under the assumptions on the initial data of Theorem  \ref{thm:sys_2}, there exists
 an interval $\jin:=[0,t_\star]$ on which the system \eqref{sys_limite:1_bis} -- \eqref{sys_limite:4_bis} admits a weak
solution $(\rho, \, f, \,\uu)$ (cf. Remark \ref{we} for a definition of weak solutions).
Furthermore $\rho$ and $f$ are nonnegative, and
\begin{align*}
\rho &\in \L^{\infty}(\jin;\L^{5/3}(\Tt)),\\
f &\in \L^{\infty}(\jin \times \Tt \times \Rt),\\
\uu &\in \L^{2}(\jin;\Hh^{1}(\Tt)) \cap \L^{\infty}(\jin;\L^{2}(\Tt)).
\end{align*}
\end{prop}
\medskip

{\bf{Proof of Proposition \ref{prop:local} }}:

 We begin with the Proof of bounds for the
solution of \eqref{eq:sysreg1} -- \eqref{eq:sysreg4}, which do not depend on $\ep$.
\smallskip

In the sequel $\Pp_T(\uuin,\rhoin,\fin)$ [or more simply $\ptin$] will denote a polynomial function (with positive coefficients) of the quantity $\|\sqrt{1+\rhoin}\,\uuin\|_2 + \|\fin\|_\infty + M_2 \fin$. Note
that $\ptin$ will always be independent of $\ep$, but may depend on $T$.

\subsection{Uniform bounds with respect to $\ep$}\label{subsec:bounds}

Let us first recall a classical Lemma linking velocity moments, the proof of which may be found in \cite{article:Ayman} (Lemma~1, Section 3.3):

\begin{Lem}\label{lemme:moments}
 Let $\gamma >0$ and $h$ be a nonnegative element of $ \textnormal{L}^\infty([0, T] \times \T^3\times \R^3)$, such that $m_\gamma h(t, \xx)< +\infty$ for a.e. $(t,\xx)$. The following estimate holds for $0\leq \alpha <\gamma$:
\begin{eqnarray*}
m_\alpha h(t,\xx)\leq \left(\frac{4}{3} \pi \|h(t, \xx, \cdot)\|_{\textnormal{L}^{\infty}(\R^3_{\xxi})}+1 \right)m_\gamma h(t,\xx)^{\frac{\alpha+3}{\gamma+3}}.
\end{eqnarray*}
\end{Lem}

Thanks to the maximum principle in eq. \eqref{eq:sysreg3}, we observe that
\begin{align}
\label{ineq:boundf}\| \fep\|_\infty &\leq \|\fin\|_{\infty}\, e^{2T},
\end{align}
from which, thanks to Lemma \ref{lemme:moments}, we deduce
\begin{align}   
\label{eq:bornem0f}\norm{m_0 \fep(t)}_{5/3} &\leq \ptin  \, M_2 \fep (t)^{3/5},\\ 
\label{eq:bornem1f}\norm{\mm_1 \fep(t)}_{5/4} &\leq \ptin \, M_2 \fep (t)^{4/5}. 
\end{align}  
Since $f_\ep$ is compactly supported and solves eq. \eqref{eq:sysreg3} strongly, one gets by
 multiplying this equation by $|\xxi|^2$ and integrating in $\xx,\xxi$: 
\begin{align}
\frac{\dd}{\dd t } M_2 \fep + 3 M_2 \fep = 2 \int_{\T^3}  (\uep \star \fiep) \cdot \mm_1(\fep)\, \dd \xx.\label{ineq:nrjfep2} 
\end{align}
Using the previous estimates, we hence have
\begin{align}
\frac{\dd}{\dd t } M_2 \fep + 3 M_2 \fep &\leq   \norm{\uep(t)}_{5} \, \stackrel{=1}{\overbrace{\norm{\fiep}_{1}}} \, \norm{\mm_1 \fep(t)}_{5/4} \nonumber \\
 & \leq  \ptin \,\norm{\uep(t)}_{5}\,  M_2 \fep(t)^{4/5}, \nonumber
\end{align}
from which we easily deduce 
\beno 
M_2 \fep (t) \leq \left\{ M_2 \fin^{\ep}  + \ptin \intt{0}{t} \norm{\uep(s)}_{5} \dd s \right\}^5.
\eeno
Using the Sobolev injection $\Hf^1(\T^3)\hookrightarrow \Lf^5(\T^3)$, one gets 
\ben \label{eq:borne_M2f}
M_2 \fep (t) \leq \ptin \left\{ 1+ \intt{0}{t} \norm{\uep(s)}_2 \dd s  + \intt{0}{t} \norm{\na_{\xx} \uep(s)}_{2} \dd s  \right\}^5.
\een
Using the previous estimate with \eqref{eq:bornem0f} -- \eqref{eq:bornem1f}, we obtain
\begin{align}   
\label{eq:bornem0f2}\norm{m_0 \fep(t)}_{5/3} &\leq\ptin \left\{ 1+ \intt{0}{t} \norm{\uep(s)}_2 \dd s  + \intt{0}{t} \norm{\na_{\xx} \uep(s)}_{2} \dd s  \right\}^3,\\ 
\label{eq:bornem1f2}\norm{\mm_1 \fep(t)}_{5/4} &\leq \ptin \left\{ 1+ \intt{0}{t} \norm{\uep(s)}_2 \dd s  + \intt{0}{t} \norm{\na_{\xx} \uep(s)}_{2} \dd s  \right\}^4. 
\end{align}  
The energy estimate \eqref{ineq:nrjep}  is satisfied with $\nu:=\rhoin^\ep$, $\oom:=\uuin^\ep$, and rewrites 
\begin{align*}
\frac{1}{2} \left\{\int_{\T^3}(1+\rho_\ep(t,\xx))|\uu_\ep(t,\xx)|^2\dd \xx\right\} &+ \int_0^t\int_{\T^3}|\nabla_\xx \uu_\ep(s,\xx)|^2\,\dd \xx\,\dd s \\
&\leq-\frac{3}{2}\int_0^t\int_{\T^3}|\uu_\ep(s,\xx)|^2 m_0(f_\ep\gamma)(s,\xx)\,\dd \xx\,\dd s\\
&+2\int_0^t\int_{\T^3}\uu_\ep(s,\xx)\cdot \mm_1(f_\ep\gamma)(s,\xx)\,\dd \xx\,\dd s\\
&+\frac{1}{2} \left\{\int_{\T^3}(1+\rhoin^\ep(\xx))\,|\uuin^\ep(\xx)|^2\,\dd \xx\right\},
\end{align*}
so that (since $\rho_\ep$ is nonnegative) we obtain
\begin{align*}
\frac{1}{2}\,\|\uu_\ep(t)\|_{2}^2 + \int_0^t \|\nabla_\xx \uu_\ep(s)\|_2^2\,\dd s & \leq  \frac{1}{2}\|\sqrt{1+\rhoin^\ep}\, \uuin^\ep\|_{2}^2 + 2 \int_0^t \|\uu_\ep\cdot \mm_1 f_\ep(s)\|_1 \,\dd s \\
&\leq  \ptin + 2 \int_0^t \|\uu_\ep\cdot \mm_1 f_\ep(s)\|_1\, \dd s. \\
\end{align*}
Using H\"older and Young inequalities, we get
\begin{align*}
\frac{1}{2}\|\uu_\ep(t)\|_{2}^2 + \frac{1}{2}\int_0^t \|\nabla_\xx \uu_\ep(s)\|_2^2\,\dd s 
&\leq \ptin + \ptin \int_0^t \|\uu_\ep(s)\|_2^2\, \dd s+  \ptin \int_0^t \|\mm_1 f_\ep(s)\|_{5/4}^2 \,\dd s.
\end{align*}
We then obtain thanks to \eqref{eq:bornem1f2}
\begin{align*}
\hspace{-1cm}\|\uu_\ep(t)\|_{2}^2 + \int_0^t \|\nabla_\xx \uu_\ep(s)\|_2^2\, \dd s
&\leq \ptin+ \ptin\int_0^t\left\{\|\uu_\ep(s)\|_2^2 + \intt{0}{s} \norm{\na \uep(\sigma)}_{2}^2\, \dd \sigma  \right\}^4\dd s.
\end{align*}
Denote by $\tilde{z}$ the maximal solution of the Cauchy problem 
$$ z'(t) = \ptin\, z(t)^4, \qquad z(0)=\ptin, $$
defined on some maximal interval $\iin$.
 Then if $\jin$ is the closure of $\iin/2$, using the nonlinear Gr\"onwall Lemma~\ref{lem:nlgron}, we get for $t\in \jin$:
\begin{align*}
\norm{\uep(t)}_2^2  +  \intt{0}{t}  \norm{\na \uep(s)}_{2}^2 \dd s\leq \|\tilde{z}\|_{\L^\infty}:=\ktin.
\end{align*}
\begin{rmk}\label{rmk:t_star}
Notice that thanks to Lemma \ref{lem:vie} (changing $\ptin$ if necessary) we always have $|\jin| \geq 1/\ptin$.
\end{rmk}

We deduce from the previous local estimate that $(\uu_\ep)_\ep$ is bounded in  $\L^\infty\big(\jin;\L^2(\T^3)\big)\cap \L^2\big(\jin;\Hh^1(\T^3)\big)$. Using \eqref{eq:borne_M2f}, we hence get the boundedness of  $(M_2\fep)_{\ep}$ in $\L^\infty(\jin)$ and then,
 with \eqref{eq:bornem0f} and \eqref{eq:bornem1f}, we see that $(m_0 \fep)_{\ep}$ and $(\mm_1 \fep)_{\ep}$ are respectively bounded in $\L^{\infty}\big(\jin;\L^{5/3}(\Tt)\big)$ and $ \L^{\infty}\big(\jin;\L^{5/4}(\Tt)\big)$. 
Thanks to a classical transport estimate, using \eqref{eq:sysreg3}, we see that $(\rhoep)_{\ep}$ is bounded in  $\L^{\infty}\big(\jin;\L^{5/3}(\Tt)\big)$.
 From \eqref{eq:borne_M2f}, we deduce that for all $\ep>0$ and $t\in \jin$, 
\begin{align}
\label{eq:borne_M2fep}
M_2{\fep}(t) \leq  \ktin.
\end{align}
 
\subsection{Compactness properties}

Recalling Proposition \ref{regul} of existence of a solution to the regularized problem (for a given $\ep$), we deduce
from the previous (uniform in $\ep$) bounds 
the existence of $\rho\in\L^\infty(\jin;\L^{5/3}(\T^3))$, $f\in\L^\infty(\jin;\L^\infty(\T^3\times\R^3))$ such that,
 up to a subsequence,
\begin{align*}
(\rho_\ep)_\ep\ &\convw{\ep}{0} \rho\text{ in }\L^\infty(\jin;\L^{5/3}(\T^3))-\star, \\
(f_\ep)_\ep\ &\convw{\ep}{0} f\text{ in }\L^\infty(\jin;\L^{\infty}(\T^3\times\R^3))-\star. \\
\end{align*}
 Using the boundedness of $(M_2 f_\ep)_\ep$ in $\L^\infty(\jin)$ and adding a subsequence if necessary, one manages to also show that
\begin{align*}
(m_0 f_\ep)_\ep \text{ and }(m_0 (f_\ep\gamma_\ep))_\ep &\convw{\ep}{0} m_0 f\text{ in }\L^\infty(\jin;\L^{5/3}(\T^3))-\star, \\
(\mm_1 f_\ep)_\ep\text{ and }(\mm_1 (f_\ep\gamma_\ep))_\ep &\convw{\ep}{0} \mm_1 f\text{ in }\L^\infty(\jin;\Lf^{5/4}(\T^3))-\star.
\end{align*}
Finally, using the bounds on $(\uu_\ep)_\ep$ and the Sobolev injection $\Hf^1(\T^3) \hookrightarrow \Lf^6(\T^3)$, we get
 the existence of $\uu\in\L^\infty(\jin;\Lf^2(\T^3))\cap\L^2(\jin;\Hf_{\div}^1(\T^3))$ such that
\begin{align*}
(\uep)_\ep \convw{\ep}{0}  \uu &\text{ in } \L^2\big(\jin;\Hf^1(\T^3)\big),\\
(\uep)_\ep \convw{\ep}{0}  \uu &\text{ in } \L^2\big(\jin; \Lf^6(\T^3)\big),\\
(\uep)_\ep  \convw{\ep}{0}  \uu &\text{ in } \L^\infty\big(\jin;\Lf^2(\T^3)\big)-\star,
\end{align*}
and of course $(\uu_\ep\star\ffi_\ep)_\ep$ converges also towards $\uu$ for the same topologies.

\subsection{Weak convergence of nonlinear terms}

In view of the previous weak convergences, it only remains to check that the nonlinear terms
 converge to the expected limits, so that the weak limit triplet $(f,\rho,\uu)$ will indeed be a solution of our system. 
Notice that since all the bounds and weak extractions are performed on the local interval $\jin$, the constructed solution will only be local in time. In the last subsection, we shall explain how to extend it. Until then, through all the current subsection, the index $_t$ will refer to the intervall $\jin$ in the notation ``$\L_t^p(\textnormal{E})$''.

\vspace{2mm}

Notice that the system can be written 
\begin{align}
\label{eq:rap0} & \div_\xx  \uu_\ep = 0,\\
\label{eq:rap1} & \partial_t f_\ep + \textnormal{div}_{\xx,\xxi}(\aaa_\ep f_\ep) - 2\,f_\ep = 0,\\
\label{eq:rap2} & \partial_t \rho_\ep + \div_\xx(\rho_\ep(\uu_\ep\star\ffi_\ep)) = m_0(f_\ep\gamma_\ep),\\
\label{eq:rap3} &\mathbb{P}\Big\{\partial_t[(1+\rho_\ep)\uu_\ep] + \div_\xx \big[(1+\rho_\ep)(\uu_\ep\star\ffi_\ep) \otimes \uu_\ep\big]\Big\} -\Delta_\xx \uu_\ep = \P\Big\{2\mm_1(f_\ep\gamma_\ep)-\uu_\ep m_0(f_\ep\gamma_\ep)\Big\},
\end{align}
where $\aaa_\ep(t,\xx,\xxi):=(\xxi,[\uu_\ep(t)\star \ffi_\ep](\xx))$. \vspace{2mm}\\
In order to handle the nonlinear terms, we shall use several times the Proposition~\ref{prop:annex_compensation} of the Appendix. 
\smallskip

We denote by $\mathscr{M}_s$ the vector space of bounded measures on $\jin\times\T^3$. 

\vspace{2mm}

\textsf{\emph{All the coming facts and their proofs are true up to some (finite number of) extractions that we don't mention in the sequel.}}

\vspace{2mm}

\begin{itemize}

\item[\textsf{\textbf{Fact 1}} : ] \emph{The products $(\rho_\ep\,(\uu_\ep\star\ffi_\ep))_\ep$ and  $(\rho_\ep \,\uu_\ep)_\ep$ both tend to $\rho\,\uu$ in $\mathscrbf{M}_s$.}

\vspace{2mm}

{\sl{Proof of Fact~1}}: Using Proposition \ref{prop:annex_compensation} of the Appendix, these two terms are handled in the same way: the velocity term (convoluted by $\ffi_\ep$ or not)  plays the role of $(a_\ep)_\ep$, and is bounded in $\L^2_t(\W^{1,2}_\xx)$, whereas $(\rho_\ep)_\ep$ plays the role of $(b_\ep)_\ep$ and is bounded in $\L^\infty_t(\L^{5/3}_\xx)$ ($(5/3)'=5/2<2^\star=6$). Thanks to \eqref{eq:rap2} and the previous bounds (see subsection \ref{subsec:bounds}), $(\partial_t \rho_\ep)_\ep$ is bounded in $\L^2_t(\Hh^{-m}_\xx)$, (where $m$ is taken large enough). We hence have
\begin{align*}
(\rho_\ep\,(\uu_\ep\star\ffi_\ep))_\ep&\operatorname*{\rightharpoonup}_{\ep\rightarrow 0} \rho\, \uu,\text{ in }\mathscrbf{M}_s-\text{w}\star,\\
(\rho_\ep \,\uu_\ep)_\ep &\operatorname*{\rightharpoonup}_{\ep\rightarrow 0} \rho\, \uu,\text{ in }\mathscr{M}_s-\text{w}\star.
\end{align*}
\item[\textsf{\textbf{Fact 2}} : ] \emph{$\displaystyle \big\langle (\rho_\ep+1)\uu_\ep,\uu_\ep\big\rangle_{\Lf^2_{t,\xx}} \operatorname*{\longrightarrow}_{\ep\rightarrow 0} \big\langle (\rho+1)\uu ,\uu \big\rangle_{\Lf^2_{t,\xx}}$.}

\vspace{2mm}

{\sl{Proof of Fact~2}}: We wish to prove that
\begin{align}
\label{conv:rhou} \operatorname*{\lim}_{\ep \rightarrow 0}  \int_{\jin} \int_{\T^3} (\rho_\ep+1)\, |\uu_\ep|^2 \,\dd \xx \,\dd t= \int_{\jin} \int_{\T^3} (\rho+1)\, |\uu|^2\, \dd \xx \,\dd t.
 \end{align}
First write, since $\uu_\ep$ is divergence-free,
\begin{align*}
\int_{\jin} \int_{\T^3} (\rho_\ep+1) \,|\uu_\ep|^2 \, \dd \xx\, \dd t &= \int_{\jin} \int_{\T^3} (\rho_\ep+1)
\, \uu_\ep \cdot \uu_\ep \, \dd \xx\, \dd t \\
&= \int_{\jin} \int_{\T^3} \mathbb{P}\big[(\rho_\ep+1) \uu_\ep\big] \cdot \uu_\ep\, \dd \xx\, \dd t.
\end{align*}
Harmonic analysis and  singular integral theory allow to show that $\P$ is bounded from $\L^q_t(\L^p_{\xx})$ to itself, for all $p\in]1,\infty[$ and $q\in[1,\infty]$ (see for instance \cite{stein}).  Since the strong continuity of an operator implies its  weak sequential continuity, we see that $\P$ is sequentially continuous from $\L^2_t(\L^{30/23}_\xx)$ to itself, equipped with the weak topology. But we have seen in Subsection \ref{subsec:bounds} that $(\uu_\ep)_\ep$ is bounded in $\L^2_t(\Hf^1_{\xx})$, which is embedded in $\L^2_t(\Lf^6_{\xx})$ by Sobolev injection, and $(\rho_\ep)_\ep$ is bounded in $\L^\infty_t(\L^{5/3}_{\xx})$ so that by H\"older inequality, $((\rho_\ep+1)\uu_\ep)_\ep$ is bounded in $\L^2_t(\Lf^{30/23}_\xx)$ and hence (up to a subsequence) $\Big(\mathbb{P}\big[(\rho_\ep+1)\, \uu_\ep\big]\Big)_\ep$ converges weakly to $\mathbb{P}\big[(\rho+1)\, \uu\big]$ in $\L^2_t(\L^{30/23}_\xx)$. We now can use Proposition \ref{prop:annex_compensationbis} of the Appendix, with $q=2$ and $r=\infty$. Indeed,
\vspace{2mm}

\begin{itemize}
\item $(\uu_\ep)_\ep$ is bounded in $\L^2_t(\W^{1,2}_\xx)\cap\L^\infty_t(\L^2_\xx)$,
\item since $\partial_t \P = \P \partial_t$, \eqref{eq:rap3} implies that $\Big(\partial_t \mathbb{P}\big[(\rho_\ep+1) \uu_\ep\big]\Big)_\ep$ is bounded in $\L^{1}_t(\Hh_\xx^{-m})$ for $m$ large enough,
\item  $(30/23)'=30/7<2^\star=6$.
\end{itemize}
\vspace{2mm}
Hence the product $\Big\{\mathbb{P}\big[(\rho_\ep+1) \,\uu_\ep\big] \cdot \uu_\ep\Big\}_\ep$ converges to $\mathbb{P}\big[(\rho+1) \,\uu\big] \cdot \uu$ in $\mathscrbf{M}_s-\star$. In particular, using
 $\mathds{1}_{\jin\times\T^3}$ as a test function, we get Fact 2.
\vspace{2mm}
\item[\textbf{\textsf{Fact 3}}] \emph{$(\uu_\ep)_\ep$ and $(\uu_\ep\star\ffi_\ep)_\ep$  both strongly converge to $\uu$ in $\Lf^2_{t,\xx}$.}

\vspace{2mm}

{\sl{Proof of Fact~3}}: We estimate
\begin{align*}
\int_{\jin} \int_{\T^3}|\uu_\ep-\uu|^2\, \dd \xx\, \dd t &\leq \int_{\jin} \int_{\T^3}(1+\rho_\ep)\,|\uu_\ep-\uu|^2 \,\dd \xx\, \dd t\\
& =  \int_{\jin} \int_{\T^3}(1+\rho_\ep)\,|\uu_\ep|^2 \,\dd \xx \,\dd t + 
 \int_{\jin} \int_{\T^3}(1+\rho_\ep)\,|\uu|^2 \,\dd \xx\, \dd t \\
&-  2\int_{\jin} \int_{\T^3}(1+\rho_\ep)\,\uu_\ep\cdot \uu \, \dd \xx \,\dd t. \end{align*}
The first term of the second line converges to the same expression, but without $\ep$: this
 is exactly Fact 2 proven above. We have the same behavior for the second term of this line:
 $(\rho_\ep)_\ep$ converges weakly to $\rho$ in $\L^\infty_t(\L^{5/3}_\xx)-\star$ and $|\uu|^2 \in \L^1_t(\L^3_\xx)$. Eventually, for the third and last term, we use the (already proven) convergence of $[(1+\rho_\ep)\,\uu_\ep]_\ep$ towards $(1+\rho)\,\uu$ in $\L^2_t(\Lf^{30/23}_\xx)$,
 and the embedding $\uu\in\L^2_t(\Lf^6_\xx)\hookrightarrow\L^2_t(\Lf^{30/7}_{\xx})$.
 Strong convergence of $(\uu_\ep)_\ep$ is then directly transfered to $(\uu_\ep\star\ffi_\ep)_\ep$.
\vspace{2mm}
\item[\textbf{\textsf{Fact 4}}] \emph{$(\uu_\ep)_\ep$ converges strongly to $\uu$ in all $\L^{c}_t(\Lf^{d}_\xx)$, for $c<2$ and $d<6$.}

\vspace{2mm}

{\sl{Proof of Fact~4}}: Thanks to Fact 3, we have strong convergence of $(\uu_\ep)_\ep$,
 and hence almost everywhere convergence. Since this family is bounded $\L^2_t(\Lf^6_\xx)$,
 we get the desired convergences.
\vspace{2mm}

\item[\textbf{\textsf{Fact 5}}] \emph{$(f_\ep\,(\uu_\ep\star\ffi_\ep))_\ep$ and 
$(\uep \, m_0(\fep \gaep))_\ep$ converges weakly in $\Lf^1$ to respectively $f\, \uu$  and $\uu \,m_0 (f)$.}

\vspace{2mm}

{\sl{Proof of Fact~5}}: The family $(f_\ep)_\ep$ converges weakly to $f$ in $\L^\infty_{t}(\L^\infty_{\xx,\xxi})-\star$, which with the strong convergence (Fact 3) of $(\uu_\ep\star\ffi_\ep)_\ep$ in $\Lf^2_{t,\xx}\hookrightarrow \Lf^1_{t,\xx}$, ensures the weak $\Lf^1_{t,\xx,\xxi}$ convergence. Similarly, $(m_0(\fep\,\gaep))_\ep$ converges weakly to $m_0 f$ in $\L^\infty_t(\L^{5/3}_\xx)-\star$, and Fact 4 ensures for instance that $(\uep)_\ep$ converges strongly in $\L^1_t(\Lf^{5/2}_\xx)$.

\vspace{2mm}

\item[\textbf{\textsf{Fact 6}}] \emph{$\big[(1+\rho_\ep)\, (\uu_\ep\star\ffi_\ep)\otimes\uu_\ep\big]_\ep$ converges weakly in $\Lf^1_{t,\xx}$ to $(1+\rho)\,\uu\otimes\uu$.}

\vspace{2mm}

{\sl{Proof of Fact~6}}: Once again, it's a `` weak $\times$ strong '' type of convergence here. First notice
\begin{align*}
\frac{2}{11} = \frac{\theta}{6} + \frac{1-\theta}{2},
\end{align*}
where $\theta:=21/22$. Hence, by H\"older's inequality, we get the following interpolation result: $\L^\infty_t(\L^2_\xx)\cap\L^2_t(\L^6_\xx)\hookrightarrow \L^{2/\theta}_t(\L^{11/2}_\xx)$. This implies that both $(\uu_\ep)_\ep$ and $(\uu_\ep\star\ffi_\ep)_\ep$ are bounded in $\L^{2/\theta}_t(\L^{11/2}_\xx)$, and the product $[(\uu_\ep\star\ffi_\ep) \otimes \uu_\ep]_\ep$ (which converges almost everywhere to $\uu\otimes\uu$, see Fact 4) is hence bounded in $\L^{1/\theta}_t(\L^{11/4}_\xx)$. Since $\theta <1$ and $11/4>5/2$, 
we get the strong convergence of this product to $\uu\otimes\uu$, in $\L^1_t(\L^{5/2}_\xx)$, and we already know that $(\rho_\ep)_\ep$ converges weakly to $\rho$ in $\L^{\infty}_t(\L^{5/3}_\xx)-\star$.
\end{itemize} 
\medskip

Using all the weak convergences above, we end up with a local in time solution of our system, which
concludes the Proof of Proposition \ref{prop:local}.
$\square$.
\medskip

Next subsection is devoted to the prolongation of the local solutions defined above in solutions
defined on $[0,T]$. 

\subsection{Energy estimate and global existence}

In order to prove the existence of global solutions to our system, 
a standard strategy consists in reproducing the previous step at time $t_\star-\ep$, obtain another local solution and paste it with the previous one, and so on and so forth.
\par
 For such a strategy to succeed, one must ensure that the sequence of local times of existence does not decrease too quickly. But, as noticed in Remark \ref{rmk:t_star}, the local time of existence is bounded
below  by $1/\ptin$, which is (by definition) a non-increasing function of $\|\sqrt{1+\rhoin}\uuin\|_2 + \|\fin\|_\infty + M_2 \fin$.
 We already know (by maximum principle) that $\|f\|_\infty \leq e^{2t_\star}\fin$, and this bound ensures
 that all possible extension will always satisfy $\|f\|_\infty \leq e^{2T}\fin$, since the sum of all local times of existence does not exceed $T$ (well if it does, we're done !). It is hence sufficient to prove that, for almost all $t\in\jin$, 
\begin{align}
\label{ineq:dec}M_2 f(t) + \|\sqrt{1+\rho(t)}\,\uu(t)\|_2 \leq M_2 \fin + \|\sqrt{1+\rhoin}\,\uuin\|_2.
\end{align}
Indeed, such an estimate would propagate for each local solution and we may hence bound from below
 all the corresponding times of existence, which means that our strategy would end in a finite number of steps. In fact \eqref{ineq:dec} is a straightforward consequence of the following energy estimate:

\vspace{2mm}
\begin{prop}\label{propo:energy}
The solution built in Proposition \ref{prop:local} satisfies for almost all $t\in\jin$,
\begin{align}
\frac{1}{2}\left\{M_2 f (t) + \| \sqrt{1+\rho(t)} \, \uu(t)\|_{\textnormal{L}^2(\T^3)}^2\right\} &+ \int_0^t\|\nabla_\xx \uu(s)\|_{\textnormal{L}^2(\T^3)}^2 \dd s + \frac{3}{2}\int_0^t\int_{\T^3 \times \R^3}|\uu-\xxi|^2 \,f\, \dd \xxi\, \dd\xx\, \dd s  \nonumber \\
& \leq \frac{1}{2}\left\{M_2 f_{\text{\bf in}} + \| \sqrt{1+\rho_{\text{\bf in}}} \,\uu_{\text{\bf in}}\|_{\textnormal{L}^2(\T^3)}^2 \right\}. \label{eqn:energie}
\end{align}
\end{prop}

{\bf{Proof of Proposition \ref{propo:energy}}}:
\medskip

We first use inequality \eqref{ineq:nrjep} for solutions of the regularized system,
 and add the integral in time of \eqref{ineq:nrjfep2}$\times \frac12$ to get
\begin{align*}
\frac{1}{2}\left\{M_2 f_\ep (t) + \| \sqrt{1+\rho(t)}\, \uu_\ep(t)\|_{\textnormal{L}^2(\T^3)}^2\right\} &+ \int_0^t\|\nabla_\xx \uu_\ep(s)\|_{\textnormal{L}^2(\T^3)}^2 \dd s \\
&+ \frac{3}{2}\int_0^t\int_{\T^3 \times \R^3}|\uu_\ep-\xxi|^2\, f_\ep \,\dd \xxi\, \dd\xx\, \dd s \\
&\leq \frac{1}{2}\,\left\{M_2 f_{\text{\bf in}} + \| \sqrt{1+\rho_{\text{\bf in}}}\,  \uu_{\text{\bf in}}\|_{\textnormal{L}^2(\T^3)}^2 \right\} + R_\ep(t),
\end{align*}
where
\begin{align*}
R_\ep(t) = \frac{3}{2}\stackrel{R_\ep^1(t)}{\overbrace{\int_0^t \int_{\T^3 \times \R^3}f_\ep\,
 |\uu_\ep|^2\,(1-\gamma_\ep(\xxi))\,\dd\xxi\, \dd \xx\, \dd s}} 
&+ 2 \stackrel{R_\ep^2(t)}{\overbrace{\int_0^t \int_{\T^3 \times \R^3}f_\ep\, \xxi \cdot \uu_\ep\,(\gamma_\ep(\xxi)-1)\, \dd\xxi\, \dd \xx\, \dd s}}\\
 &+ \stackrel{R_\ep^3(t)}{\overbrace{\int_0^t \int_{\T^3 \times \R^3}f_\ep \,\xxi \cdot (\uu_\ep \star \ffi_\ep-\uu_\ep)\, \dd\xxi\, \dd \xx\, \dd s}}.
\end{align*}
Using Lemma \ref{lemme:moments}, the maximum principle, and bound \eqref{eq:borne_M2fep}, we notice that for $\alpha>0$ small enough, $(m_\alpha f_\ep)_\ep$ is bounded in $\L^\infty_t(\L^{3/2}_{\xx})$. Now recall that $\gamma_\ep$ is chosen so that $|1-\gamma_\ep|(\xxi) \leq \mathds{1}_{|\xxi|\gtrsim 1/\ep}$, hence 
\begin{align*}
R_\ep^1(t) \lesssim \ep^\alpha \|m_\alpha f_\ep\|_{\L^\infty_t(\L^{3/2}_{\xx})} \|\uu_\ep\|_{\L^2_t(\L^6_{\xx})}^2,
\end{align*}
which goes to $0$ with $\ep$ since $(\uu_\ep)_\ep$ is bounded in $\L^2_t(\Hf^1_{\xx}) \hookrightarrow \L^2_t(\Lf^6_{\xx})$. This shows that $(R_\ep^1)_\ep \ds\conv{\ep}{0} 0$ uniformly on $\jin$, and a similar Proof applies for $(R_\ep^2)_\ep$. Then, we have 
\begin{align*}
R_\ep^3(t) = \int_0^t \int_{\T^3\times\R^3} (f_\ep \xxi \cdot \vv_\ep)\, \dd \xxi\, \dd\xx\,\dd s ,
\end{align*}
where $(\vv_\ep)_\ep$ is bounded in $\L^2_t(\Lf^6_{\xx})$ and converges to $0$ in $\L^2_{t,\xx}\hookrightarrow\L^1_{t,\xx}$. For values of $\xxi$ satisfying $|\xxi|\leq\ds (\|\vv_\ep\|_1+\ep)^{-1/2}$, we simply use the maximum principle for $(f_\ep)_\ep$ to see that their contribution goes to $0$ with $\ep$ (uniformly in time), and for the large values of $\xxi$ satisfying the opposite inequality, we handle them as we did for $R_\ep^1$ and $R_\ep^2$, using the $\L^2_t(\Lf^6_{\xx})$ bound.

\vspace{2mm}

At this stage, we proved that $(R_\ep)_\ep$ converges to $0$ uniformly in time. Let us treat the other terms of inequality \eqref{eqn:energie}. We first use the classical estimates of weak convergence to get 
\begin{align}
\nonumber\int_0^t \|\nabla_\xx \uu(s) \|_{\L^2(\T^3)}^2\dd s & \leq \operatorname*{\underline{\lim}}_{\ep \rightarrow 0} \hspace{1mm} \int_0^t \| \nabla_{\xx} \uu_{\ep} \|_{\L^2(\T^3)}^2 \dd s, \\
\label{ineq:classm2}M_2 f(t) &\leq \operatorname*{\underline{\lim}}_{\ep \rightarrow 0} \hspace{1mm} M_2 f_\ep (t).
\end{align}
Furthermore, adapting the proof of estimate \eqref{conv:rhou} used in Fact $2$,
 we obtain for almost all times $t\in\jin$,
\begin{align*}
\| \sqrt{1+\rho_\ep(t)}\,\uu_\ep(t)\|_{\L^2(\T^3)}^2 \operatorname*{\rightarrow}_{\ep \rightarrow 0} \|\sqrt{1+\rho(t)}\,\uu(t)\|_{\L^2(\T^3)}^2,
\end{align*}
so that it only remains to prove
\begin{align}
\label{preuvestilim4}\int_0^t\int_{\T^3 \times \R^3}|\uu-\xxi|^2 f \dd \xxi \dd\xx \dd s &\leq \operatorname*{\underline{\lim}}_{\ep \rightarrow 0}\int_0^t\int_{\T^3 \times \R^3}|\uu_\ep-\xxi|^2 f_\ep \dd \xxi \dd\xx \dd s,
\end{align}
in order to conclude the Proof of Proposition \ref{propo:energy}. We first write
\begin{align}
\label{preuve0} \int_0^t\int_{\T^3 \times \R^3}|\uu_\ep-\xxi|^2 f_\ep \dd \xxi \dd\xx \dd s = \int_0^t\int_{\T^3 \times \R^3}|\uu_\ep|^2 f_\ep \dd \xxi \dd\xx \dd s + \int_0^t M_2 f_\ep(s)\dd s
- 2 \int_0^t\int_{\T^3 \times \R^3} \xxi \cdot \uu_\ep f_\ep \dd \xxi \dd\xx \dd s.
\end{align}
Since $(\uu_\ep)_\ep$ converges strongly in $\Lf^2_{t,\xx}$ and $(f_\ep)_\ep$ converges weakly in $\L^\infty_{t,\xx,\xxi}-\star$, we have by Fatou's Lemma 
\begin{align*}
\int_0^t\int_{\T^3 \times \R^3}|\uu|^2 f \mathds{1}_{|\xxi|\leq  n} \dd \xxi\, \dd\xx\, \dd s \leq  \operatorname*{\underline{\lim}}_{\ep \rightarrow 0} \int_0^t\int_{\T^3 \times \R^3}|\uu_\ep|^2 f_\ep \dd \xxi \,\dd\xx \,\dd s,
\end{align*}
hence by monotone convergence,
\begin{align}
\label{preuve1}\int_0^t\int_{\T^3 \times \R^3}|\uu|^2 f \, \dd \xxi\, \dd\xx\, \dd s \leq  \operatorname*{\underline{\lim}}_{\ep \rightarrow 0} \int_0^t\int_{\T^3 \times \R^3}|\uu_\ep|^2 f_\ep \,
\dd \xxi \,\dd\xx \,\dd s. 
\end{align}
Using another time Fatou's Lemma with estimate \eqref{ineq:classm2} shows that 
\begin{align}
\label{preuve2}\int_0^t M_2 f(s) \dd s \leq  \operatorname*{\underline{\lim}}_{\ep \rightarrow 0} \int_0^t M_2 f_\ep(s) \dd s.
\end{align}
On the other hand, since $v_\ep:=\uu_\ep-\uu$ is bounded in $\L^2_t(\Lf^6_{\xx})$ and converges to $0$ in $\L^1_{t,\xx}$, we have as before (see the study of $R_\ep^3$)
\begin{align*}
\int_0^t \int_{\T^3\times\R^3} f_\ep \, \xxi \cdot (\uu_\ep-\uu)\,
\dd \xxi\,\dd\xx \,\dd s \conv{\ep}{0} 0,
\end{align*}
and handling large velocities as we did for $R_\ep^1$ and $R_\ep^2$, we can show that 
\begin{align*}
\int_{\T^3\times\R^3} (f \xxi \cdot \uu) \mathds{1}_{|\xxi| > n}\dd \xxi\,\dd\xx \,\dd s + \sup_{\ep>0}\int_0^t \int_{\T^3\times\R^3} (f_\ep \xxi \cdot \uu) \mathds{1}_{|\xxi| > n}\dd \xxi\,\dd\xx \,\dd s \conv{n}{+\infty} 0.
\end{align*}
When $n\in\N$ is fixed, we have by weak $\L^\infty_{t,\xx,\xxi}-\star$ convergence of $(f_\ep)_\ep$,  
\begin{align*}
\int_0^t \int_{\T^3\times\R^3} (f_\ep \, \xxi \cdot \uu) \mathds{1}_{|\xxi| \leq n}\, \dd \xxi\,\dd\xx \,\dd s \conv{n}{+\infty} \int_0^t \int_{\T^3\times\R^3} (f_\ep\, \xxi \cdot \uu)\, \dd \xxi\,\dd\xx \,\dd s,
\end{align*}
so that the three last convergences imply together
\begin{align}
\label{preuve3}\int_0^t \int_{\T^3\times\R^3} (f_\ep \, \xxi \cdot \uu_\ep) \, \dd \xxi\,\dd\xx \,\dd s \conv{n}{+\infty} \int_0^t \int_{\T^3\times\R^3} (f \,\xxi \cdot \uu)\, \dd \xxi\,\dd\xx \,\dd s.
\end{align}
Using \eqref{preuve1}, \eqref{preuve2} and \eqref{preuve3} in \eqref{preuve0},
 we get \eqref{preuvestilim4} and this ends the Proof of Proposition \ref{propo:energy}.
 $\square$
\medskip

This also concludes the Proof of Theorem \ref{thm:sys_2}, thanks to the strategy of prolongation of local solutions explained at the beginning of the Subsection.
$\square$

\section{Appendix} 

 We present in the Appendix a few auxiliary Lemmas which are used in various parts of this work. We start
 with a standard variant of Gr\"onwall's Lemma, that we recall for the sake of completeness.

\subsection{A variant of Gr\"onwall's Lemma} 

\begin{Lem}\label{lem:nlgron}
Let $f\in\mathscr{C}^1(\R_+)$ be a convex non-decreasing function and $\alpha\in \R$. Consider $(z,[0,t_\star[)$ the maximal solution of $z'=f(z)$, $z(0)=\alpha$ on $\R_+$. 
Let $a:\R_+\rightarrow\R_+$ be a continuous function such that
\begin{align*}
a(t) \leq \alpha + \int_0^t f(a(s))\, \dd s.
\end{align*}
Then, for all $t\in[0,t_\star[$,
\begin{align*}
a(t) \leq z(t).
\end{align*}
\end{Lem}

{\bf{Proof of Lemma \ref{lem:nlgron}}}:
\medskip

The usual (linear) Gr\"onwall lemma shows that if $u\in\mathscr{C}^1(\R_+)$ satisfies $\ds u(t) \leq \int_0^t b(s) u(s) \dd s$ on $[0,t_\star[$, then $u$ is nonpositive on this interval.
 We have here by convexity (on $[0,t_\star[$)
\begin{align*}
u(t):=a(t) -z(t) \leq \int_0^t [f(a(s))-f(z(s))]\,  \dd s \leq \int_0^t f'(a(s))\, u(s)\, \dd s,
\end{align*}
and since $f$ is nondecreasing, $u \leq 0$ on $[0,t_\star[$.
$\square$
\vspace{2mm}

As a consequence, we get the

\begin{Lem}\label{lem:vie}
Let $A,\gamma>0$. The maximal solution of the Cauchy problem $z'=A\,z^{1+\gamma}$, $z(0)=A$ is defined at least on $[0,(\gamma \,A^{\gamma+1})^{-1}[$.
\end{Lem}

{\bf{Proof of Lemma \ref{lem:vie}}}
\medskip

Indeed : $\ds t\mapsto (A^{-\gamma}- A\,\gamma\, t)^{-1/\gamma}$ is well defined for $t < (\gamma\, A^{\gamma+1})^{-1}$, and solves the Cauchy problem.
$\square$
\medskip

In next Subsection, we recall for the sake of completeness the following classical compactness Lemma:

\subsection{A compactness lemma}

\begin{Lem} \label{annex:aubin}
The injection $\Hh^1_{t,\xx}\hookrightarrow\mathscr{C}^0_t(\L^2_{\xx})$ is compact.
\end{Lem}

\vspace{2mm}

{\bf{Proof of Lemma \ref{annex:aubin}}}:

\vspace{2mm}

Let $(f_n)_n$ be a bounded sequence of $\H^1_{t,\xx}$.
\bigskip

{\it{Step 1}}:  $ A = \{t\in[0,T]\,:\, \norm{f_n(t)}_{\Hh^1_{\xx}} \rightarrow +\infty\}$ is a null set with respect to the Lebesgue measure $\mu$. Indeed, if it were not, then for any $R>0$, we would have
\begin{align*}
R\, \mu(A) \leq \int_{A} \operatorname*{\underline{\lim}}_{n \rightarrow \infty} \norm{f_n(t)}^2_{\Hh^1_\xx}(t) \, \dd t.
\end{align*}
Fatou's lemma would then imply 
\begin{align*}
R \,\mu(A) \leq \operatorname*{\underline{\lim}}_{n \rightarrow \infty} \int_{A} \norm{f_n(t)}^2_{\Hh^1_\xx}  \leq \operatorname*{\underline{\lim}}_{n \rightarrow \infty} \norm{f_n}^2_{\Hh^1_{t,\xx}},
\end{align*}
which is impossible since the right-hand side is finite. By definition of $A$, for all $t$ in $[0,T]\backslash A$, the sequence $(\norm{f_n(t)}_{\Hh^1_\xx})_n$ has a bounded subsequence. Since $\mu(A)=0$, we can find a countable subset $[0,T]\backslash A$ that is dense in $[0,T]$. Let us denote this subset $\B=(t_p)_{p\in\N}$. We can then extract a (diagonal) subsequence of $(f_n)_{n}$ (sill denoted $(f_n)_n$) such that for all $p\in\N$, the sequence $(f_n(t_p))_{n\in\N}$ is bounded in ${\Hh^1_\xx}$.  
\bigskip



{\it{Step 2}}: Since $(\partial_t f_n)_n$ is bounded in $\L^2_{t,\xx}$,  $(f_n)_n$ is  \textbf{uniformly equicontinuous} w.r.t. $t$, with values in $\L_\xx^2$.
\bigskip

{\it{Step 3}}: Since by Rellich's Theorem, the injection $\Hh^1_\xx \hookrightarrow \L^2_\xx$ is compact, we can extract (again, diagonally) a subsequence (still denoted) $(f_{n})_n$ such that for all $p$, the sequence $(f_{n}(t_p))_n$ converges in $\L^2_\xx$ to some element $f(t_p) \in \L^2_\xx$. 
\bigskip

{\it{Step 4}}: $f$ is uniformly continuous on $\B$, with values in $\L^2_\xx$. Indeed, for any $t$, $s$ in $\B$,
\begin{align*}
\norm{f(t)-f(s)}_{\L^2_\xx} \leq \norm{f(t)-f_{n}(t)}_{\L^2_\xx}+ \norm{f_{n}(t)-f_{n}(s)}_{\L^2_\xx} +\norm{f_{n}(s)-f(s)}_{\L^2_\xx}.
\end{align*}
The central term of the right-hand side goes to $0$ with $|t-s|$, independently of $n$ because of the uniform equicontinuity obtained in Step 2. When $t$ and $s$ are fixed elements of $\B$, the two other terms of the right-hand side vanish when $n\rightarrow \infty$ because of the extraction of Step 3.
\bigskip

{\it{Step 5}}:
$f$ admits a unique continuous extension to $[0,T]$, which is uniformly continuous on this interval, with values in $\L^2_\xx$. 
\bigskip

{\it{Step 6}}: $f$ is actually the limit of $(f_n)_n$ in $\C^0_t(\L^2_{\xx})$. Indeed, for $\ep>0$, choose $\delta$ a corresponding modulus of equicontinuity for $f$ and of uniform equicontinuity for $(f_n)_n$. Then pick ($\B$ is dense in the compact set $[0,T]$) a finite number $t_1,\dots,t_N$ of elements of $\B$ such 
that $[0,T] \subseteq \cup_{i=1}^N ]t_i-\delta,t_i+\delta[$. Eventually, for any $\sigma\in[0,T]$, if $\sigma\in]t_i-\delta,t_i+\delta[$, then 
\begin{align*}
\norm{f(\sigma) - f_{n}(\sigma) }_{\L^2_\xx}  &\leq \norm{f(\sigma) - f(t_i) }_{\L^2_\xx}  + \norm{f(t_i) - f_{n}(t_i) }_{\L^2_\xx} + \norm{f_{n}(t_i) - f_{n}(\sigma) }_{\L^2_\xx} \\
&\leq 2\ep   + \norm{f(t_i) - f_{n}(t_i) }_{\L^2_\xx},
\end{align*}
so that  $\norm{f(\sigma) - f_{n}(\sigma) }_{\L^2_\xx} \leq 3\ep$ for $n$ large enough. This ends the proof of Lemma \ref{annex:aubin}.
$\square$

 \subsection{Weak convergence of a product} \label{sec:annex_compensation}

We present here a result based on the method used in \cite{art:diperna_lions}.
\begin{prop}
\label{prop:annex_compensation}
Let $q \in [1, \infty]$ and $p\in [1,3[$. Consider two families $(a_\ep)_\ep$ and $(b_\ep)_\ep$ respectively in  $\L^q_t(\W^{1,p}_\xx)$ and $\L^{q'}_t(\L^{s'}_\xx)$, with $s< p^\star:=\frac{3p}{3-p}$. Assume the weak convergences
\begin{align*}
(a_\ep)_\ep&\operatorname*{\rightharpoonup}_{\ep\rightarrow 0} a,\text{ in}\quad \L^q_t(\W^{1,p}_\xx)- w \star,\\
(b_\ep)_\ep&\operatorname*{\rightharpoonup}_{\ep\rightarrow 0} b,\text{ in} \quad \L^{q'}_t(\L^{s'}_\xx)- w \star.
\end{align*}
If  $(\partial_t b_\ep)_\ep$ is bounded in $\L^{q'}_t(\Hh^{-m}_\xx)$ for some  $m\in\Z$, then, up to a subsequence, we have the weak convergence
\begin{align*}
(a_\ep \, b_\ep)_\ep \operatorname*{\rightharpoonup}_{\ep\rightarrow 0} 
a\, b\text{ in the sense of measures},
\end{align*}
\emph{i.e.} with test functions in $\mathscr{C}^0_{t,\xx}$.
\end{prop}
{\bf{Proof of Proposition \ref{prop:annex_compensation} }}:
\medskip

Let us first notice that for all $\ep>0$ $a_\ep \,b_\ep\in \L^1_{t,\xx}$ (Sobolev injection) and that the sequences $(a_\ep)_\ep$ and $(b_\ep)_\ep$ are bounded in the spaces in which they converge weakly. 

\vspace{2mm}

\begin{itemize}
\item[\emph{Step 1.}] We have
\begin{align*}
a \,(b\star\ffi_\eta) \operatorname*{\longrightarrow}_{\eta \rightarrow 0} a\, b,\text{ in }\L^1_{t,\xx} \text{ strong } .
\end{align*}
\item[\emph{Step 2.}] Since $(\partial_t b_\ep)_\ep$ is bounded in $\L^{q'}_t(\Hh^{-m}_\xx)$,  $(b_\ep\star\ffi_\eta)_\ep$ is bounded in $\W^{1,q'}_t(\W^{1,p'}_\xx)$ so that (thanks to Rellich's Theorem), for all fixed $\eta$, $(b_\ep \star\ffi_\eta)_\ep$ admits a (strongly) converging subsequence in  $\L^{q'}_t(\L^{p'}_\xx)$, the limit being necessarily $b\star\ffi_\eta$ (this is due to the 
uniqueness of the weak$-\star$ limit). In fact, we can choose (but we don't write it explicitly) a common (diagonal) extraction for all $\eta$ after discretization ($\eta:=1/k$). Since $(a_\ep)_\ep$ converges weakly in $\L^q_t(\L^p_\xx)$, we eventually get, for all fixed $\eta$,
\begin{align*}
(a_\ep \,(b_\ep\star\ffi_\eta))_\ep \operatorname*{\rightharpoonup}_{\ep \rightarrow 0} 
a\, (b\star\ffi_\eta) \text{ in } \L^1_{t,\xx}\text{ weak.}
\end{align*}
\item[\emph{Step 3.}] 
 We shall use the following ``commutator Lemma'', the Proof of which is rather close to the usual Friedrichs Lemma (which is a key element of \cite{art:diperna_lions}). 
\vspace{2mm}
\begin{Lem}\label{ll}
Under the assumptions\footnote{The assumption on $(\partial_t b_\ep)_\ep$ is obviously useless here.} of Proposition \ref{prop:annex_compensation}, if $(\ffi_\eta)_\eta$ is a sequence of even
 mollifiers, then the commutator (convolution in $\xx$ only) 
\begin{align*}
S_{\ep,\eta} &:= a_\ep\, (b_\ep\star\ffi_\eta)-(a_\ep\, b_\ep)\star\ffi_\eta
\end{align*}
goes to $0$  in $\L^1_{t,\xx}$ as $\eta\rightarrow 0$, uniformly in $\ep$.
\end{Lem}

{\bf{Proof of Lemma \ref{ll}}}:
\medskip

First recall the following standard fact : since $(a_\ep)_\ep$ is bounded in $\L^q_t(\W^{1,p}_{\xx})$ and $s<p^\star$, the sequence $(\tau_\hh a_\ep-a_\ep)_\ep$ tends to $0$ in $\L^q_t(\L^{s}_\xx)$ as $\hh\rightarrow 0$, uniformly in $\ep$. We now write the following equality for the commutator  \begin{align*}
S_{\ep,\eta}(t,\xx) =\int_{\B_\eta}\Big[a_\ep(t,\xx)-a_\ep(t,\xx-\yy)\Big]\,
b_\ep(t,\xx-\yy)\, \ffi_\eta(\yy)\, \dd \yy,
\end{align*}
whence thanks to Fubini's Theorem,
\begin{align*}
\|S_{\ep,\eta}\|_1\leq \|b_\ep\|_{\L^{q'}_t(\L^{s'}_\xx)}  \int_{\B_\eta} |\ffi_\eta(\yy)| \,
\|\tau_\yy a_\ep-a_\ep\|_{\L^{q}_t(\L^{s}_\xx)}
\, \dd \yy,
\end{align*}
which yields the desired uniform convergence, and concludes the Proof
of Lemma \ref{ll}. $\square$

\item[\emph{Step 4.}] We have 
\begin{align*}
(a_\ep b_\ep)\star\ffi_\eta -a_\ep b_\ep \operatorname*{\longrightarrow}_{\eta \rightharpoonup 0} 0 \text{ in the sense of measures,}
\end{align*}
uniformly in $\ep$ (with a fixed continuous test function). Indeed, if $\theta\in\mathscr{C}_{t,\xx}^0$, since $\ffi_\eta$ is even, we know that
\begin{align*}
\langle (a_\ep \,b_\ep)\star\ffi_\eta-a_\ep\, b_\ep,\theta\rangle = \langle a_\ep\, b_\ep,\theta\star\ffi_\eta-\theta\rangle,
\end{align*}
and the right-hand side tends to $0$ with the desired uniformity because $(a_\ep \,b_\ep)_\ep$ is bounded in $\L^1_{t,\xx}$, and $(\theta\star\ffi_\eta-\theta)_\eta$ goes to $0$ in $\L^\infty_{t,\xx}$ ($\theta$ is uniformly continuous).
\item[\emph{Step 5.}] Write
\begin{align*}
a\,b-a_\ep\, b_\ep &= a\,b- a\,(b\star\ffi_\eta)\\
&+ a\,(b\star\ffi_\eta)-a_\ep\, (b_\ep\star\ffi_\eta)\\
&+ a_\ep\,(b_\ep\star\ffi_\eta)-(a_\ep\, b_\ep)\star\ffi_\eta\\
&+ (a_\ep \,b_\ep)\star\ffi_\eta -a_\ep\, b_\ep.
\end{align*}
Fix $\theta\in\mathscr{C}^0_{t,\xx}$. In the right-hand side, line number $i\in\{1,2,3,4\}$ corresponds to the Step $i$ proven previously. We choose first $\eta$ to handle (uniformly in $\ep$) all the lines of the right-hand side, except the second one. Then, we
 choose the appropriate $\ep$ to handle the second line, thanks to Step 2.
 This concludes the Proof of Proposition \ref{prop:annex_compensation}.
$\square$ 
\end{itemize}
\vspace{2mm}

\vspace{2mm} 

\begin{prop}
\label{prop:annex_compensationbis} In Proposition \ref{prop:annex_compensation}, when
 $q<\infty$, the same conclusion holds assuming
 only a bound in $\L^1_t(\H^{-m}_{\xx})$ (instead of $\L^{q'}_t(\H^{-m}_{\xx})$) for $(\partial_t b_\ep)_\ep$ if, in addition, we assume the convergence of $(a_\ep)_\ep$ to $a$ in $\L^r_t(\L^p_{\xx})$ weak$-\star$,
 for some $r>q$.
\end{prop}

{\bf{Proof of Proposition \ref{prop:annex_compensationbis}}}:
\medskip

 The Proof is identical to the Proof of Proposition~\ref{prop:annex_compensation}, except for the second step (the only one using the bound on $(\partial_t b_\ep)_\ep$). For this step, we use a (diagonal) extraction such that,
 for all fixed $\eta$, $(b_\ep \star\ffi_\eta)_\ep$ converges strongly (and almost everywhere) to $b\star\ffi_\eta$ (but only) in $\L^{1}_t(\L^{p'}_\xx)$. Since $(b_\ep)_\ep$ (and hence $(b_\ep\star\ffi_\eta)_\ep$) is bounded in $\L^{q'}_t(\L^{p'}_\xx)$, with $q'>r'\geq 1$, we see that $(b_\ep\star\ffi_\eta)_\ep$ converges to $b\star\ffi_\eta$  strongly in $\L^{r'}_t(\L^{p'}_\xx)$, and the added assumption of weak-$\star$ convergence for $(a_\ep)_\ep$ allows hence to get 
\begin{align*}
(a_\ep \,(b_\ep\star\ffi_\eta))_\ep \operatorname*{\rightharpoonup}_{\ep \rightarrow 0} a\, (b\star\ffi_\eta) \text{ in $\L^1_{t,\xx}$ weak,}
\end{align*}
which ends the Proof of Proposition \ref{prop:annex_compensationbis}.
$\square$

\end{document}